\documentclass[lineno]{biometrika}

\usepackage{amsmath, amssymb}

\usepackage{times}
\usepackage{bm}
\usepackage{natbib}

\usepackage[plain,noend]{algorithm2e}

\makeatletter
\renewcommand{\algocf@captiontext}[2]{#1\algocf@typo. \AlCapFnt{}#2} 
\def\@algocf@capt@plain{top}
\renewcommand{\algocf@makecaption}[2]{%
  \addtolength{\hsize}{\algomargin}%
  \sbox\@tempboxa{\algocf@captiontext{#1}{#2}}%
  \ifdim\wd\@tempboxa >\hsize
    \hskip .5\algomargin%
    \parbox[t]{\hsize}{\algocf@captiontext{#1}{#2}}
  \else%
    \global\@minipagefalse%
    \hbox to\hsize{\box\@tempboxa}
  \fi%
  \addtolength{\hsize}{-\algomargin}%
}
\makeatother

\newcommand{\beq}{\begin{equation}}
\newcommand{\eeq}{\end{equation}}
\newcommand{\bey}{\begin{eqnarray}}
\newcommand{\eey}{\end{eqnarray}}
\newcommand{\argmin}{\operatornamewithlimits{argmin}}

\begin{document}

\jname{Tsao and Wu}
\jyear{{\em 2013b}}
\jvol{}
\jnum{}
\accessdate{ }
\copyrightinfo{{}}


\markboth{TSAO \and WU}{Extended empirical likelihood}

\title{Extended empirical likelihood for estimating equations}

\author{MIN TSAO \and FAN WU}
\affil{Department of Mathematics and Statistics, University of Victoria, Victoria,\\ British Columbia, Canada V8W 3R4 \email{mtsao@uvic.ca} \email{fwu@uvic.ca}}

\maketitle

\begin{abstract}
We derive an extended empirical likelihood for parameters defined by estimating equations which generalizes the original empirical likelihood for such parameters to the full parameter space. Under mild conditions, the extended empirical likelihood has all asymptotic properties of the original empirical likelihood. Its contours retain the data-driven shape of the latter. It can also attain the second order accuracy.  The first order extended empirical likelihood is easy-to-use yet it is substantially more accurate than other empirical likelihoods, including second order ones. We recommend it for practical applications of the empirical likelihood method.
\end{abstract}


\section{Introduction}

One important application of the empirical likelihood (Owen, 2001) is for inference on parameters defined by estimating equations $E[g(X;\theta)]=0$, where $g(x;\theta)\in \mathbb{R}^q$ is an estimating function for the parameter vector $\theta\in \mathbb{R}^p$ of a random vector $X\in \mathbb{R}^d$ (Qin and Lawless, 1994). The estimating equations are said to be just-determined if $q=p$ and over-determined if $q>p$. The latter case arises when extra information about the parameter is available and results in an estimating function of dimension $q>p$. In principle, extra information should increase the accuracy of the inference. However, Qin and Lawless (1994) noted that empirical likelihood confidence regions for over-determined cases can have substantial undercoverage.

The poor accuracy of empirical likelihood confidence regions have also been noted by others, {\em e.g.,} Hall and La Scala (1990), Corcoran, Davison and Spady (1995), Owen (2001), Tsao (2004) and Chen, Variyath and Abraham (2008). In particular, Corcoran, Davison and Spady (1995) observed that higher-order empirical likelihood method also performs poorly for small and moderate samples, suggesting that the underlying cause of the poor accuracy is not the asymptotic order of the method.  The main culprit turns out to be the mismatch between the domain of the empirical likelihood and the parameter space (Tsao, 2013; Tsao and Wu, 2013); whereas the parameter space is in general the entire $\mathbb{R}^p$, the domain is usually a bounded subset of $\mathbb{R}^p$. This mismatch is a consequence of a convex hull constraint embedded in the formulation of the empirical likelihood; values of $\theta\in \mathbb{R}^p$ that violate this constraint are excluded from the domain, leading to the mismatch. There are three variations of the original empirical likelihood (OEL) of Owen (1990) that tackle the convex hull constraint in different ways: [1] the penalized empirical likelihood (PEL) of Bartolucci (2007) and Lahiri and Mukhopadhyay (2012), [2] the adjusted empirical likelihood (AEL) by Chen, Variyath and Abraham (2008), Emerson and Owen (2009), Liu and Chen, (2010) and Chen and Huang (2012), and [3] the extended empirical likelihood (EEL) of Tsao (2013) and Tsao and Wu (2013). The PEL replaces the convex hull constraint in the OEL with a penalizing term based on the Mahalanobis distance. The AEL adds one or two pseudo-observations to the sample to ensure the convex hull constraint is never violated. The EEL expands the OEL domain geometrically to overcome the constraint and the mismatch. The AEL is available for parameters defined by estimating equations. The PEL and EEL on $\mathbb{R}^p$ are only available for the mean. The AEL, PEL and EEL all have the same asymptotic distribution as the OEL, but the EEL is a more natural generalization of the OEL as it also has identically shaped contours as the OEL. The data-driven shape of the OLE contours is a celebrated advantage of the empirical likelihood method. The EEL retains this important advantage.

In this paper, we generalize the results of Tsao and Wu (2013) for the mean to derive an EEL on $\mathbb{R}^p$ for the large collection of parameters defined by estimating equations. Under certain conditions, this EEL also has the same asymptotic properties and identically shaped contours as the OEL. It can also attain the second order accuracy of the Bartlett corrected empirical likelihood (BEL) of DiCiccio, Hall and Romano (1991). We highlight the first order version of this EEL which is not only easy-to-use but also substantially more accurate than the OEL. Surprisingly, it is also more accurate than available second order empirical likelihood methods. Because of its simplicity and accuracy, we recommend it to practitioners of the empirical likelihood method. Apart from obtaining the EEL on $\mathbb{R}^p$ for the large collection of parameters defined by estimating equations, a secondary objective of this paper is to provide through the supplementary material details of techniques for deriving the EEL on $\mathbb{R}^p$ which may be applied to parameters beyond the standard estimating equations framework. For brevity, we will use ``OEL $l(\theta)$'' and ``EEL $l^*(\theta)$'' to refer to the original and extended empirical log-likelihood ratios, respectively.




\section{Extended empirical likelihood for estimating equations}

\subsection{Preliminaries}

Let $X\in\mathbb{R}^d$ be a random vector with a parameter $\theta_0\in \mathbb{R}^p$. Let $g(X,\theta)$ be a $q$-dimensional estimating function for $\theta_0$ satisfying $E[g(X,\theta_0)]=0$ and let $X_1, \dots,$ $X_n$ be $n$ independent copies of $X$ where $n>q$. 
For simplicity, we assume all three conditions below hold.  \vspace{0.1in}

{\em Condition 1.} $E[g(X,\theta_0)]=0$ and $V[g(X,\theta_0)]\in \mathbb{R}^{q\times q}$ is positive definite.

{\em Condition 2.} $\partial g(X,\theta)/\partial \theta$ and $\partial g^2(X,\theta)/\partial \theta \partial \theta^T$ are both continuous in $\theta$, and for $\theta$ in a neighbourhood of $\theta_0$, they are both bounded in norm by some integrable function of $X$. 

{\em Condition 3.} $\lim \sup_{\|t\|\rightarrow \infty}|E[\exp\{it^Tg(X,\theta_0)\}]| <1$ and $E\|g(X,\theta_0)\|^{15}<\infty$.\\[0.1in]
These ensure the OEL for estimating equations is Bartlett correctable.
See Chen and Cui (2007) and Liu and Chen (2010).
The original empirical likelihood ratio for a $\theta \in \mathbb{R}^p$ is 
\beq
R(\theta) = \sup \left\{ \prod_{i=1}^{n} nw_i | 
\sum_{i=1}^nw_ig(X_i,\theta)=0, w_i \geq 0, \sum_{i=1}^nw_i=1 \right\}, \label{2.10}
\eeq
where $0$ is the origin in $\mathbb{R}^q$. See Owen (2001) and Qin and Lawless (1994). The OEL $l(\theta)$ is given by $l(\theta)=-2\log R(\theta)$. 
Denote by $\bar{w}=(w_1,\dots,w_n)$ a weight vector with strictly positive weights where $w_i>0$ and $\sum_{i=1}^n w_i=1$.
The domain $\Theta_n$ of the OEL $l(\theta)$ is given by
\beq
\Theta_n=\{\theta: \mbox{ $\theta\in \mathbb{R}^p$ and there exists $\bar{w}$ such that $\sum_{i=1}^n w_ig(X_i,\theta)=0$}\}. \label{2.11}
\eeq
We assume without loss of generality that $\Theta_n$ is a non-empty open set in $\mathbb{R}^p$ (see Appendix).

For a $\theta \in \Theta_n$, applying the method of Lagrange multipliers, we have
\beq
l(\theta) = 2 \sum^n_{i=1} \log\{1+\lambda^Tg(X_i,\theta)\}, \label{2.15}
\eeq
where the multiplier $\lambda=\lambda(\theta)\in \mathbb{R}^q$ satisfies
\beq
\sum_{i=1}^n \frac{g(X_i,\theta)}{1+\lambda^Tg(X_i,\theta)} = 0. \label{2.20}
\eeq
Owen (1990, 2001) showed that
$l(\theta_0)$ converges in distribution to a $\chi^2_q$ random variable as $n$ goes to infinity. Thus, the 100($1-\alpha$)\% OEL confidence region for $\theta_0$ is
\beq
{\mathcal C}_{1-\alpha}=\{\theta: \theta \in \Theta_n \mbox{ and } l(\theta) \leq c\},   \label{2.30}
\eeq
where $c$ is ($1-\alpha$)th quantile of the $\chi^2_q$ distribution. The coverage error of ${\mathcal C}_{1-\alpha}$ is given by
\beq
pr(\theta_0\in{\mathcal C}_{1-\alpha})= pr[l(\theta_0)\leq c]=pr(\chi^2_q \leq c)+O(n^{-1}). \label{2.35}
\eeq
We now briefly review the Bartlett correction (DiCiccio, Hall and Romano, 1991) for $l(\theta)$.
Under the three conditions, it can be shown that $l(\theta_0)$ has the following expansion
\beq
l(\theta_0)=nR^TR+O_p(n^{-3/2}),\label{2.36}
\eeq 
where $R$ is a $q$-dimensional vector which is a smooth function of general means. 
Through an Edgeworth expansion for the density of $n^{1/2}R$, we can show
\beq
pr\{nR^TR[1-bn^{-1}+O_p(n^{-3/2})]\leq c\}=pr(\chi^2_q \leq c)+O(n^{-2}), \label{2.37}
\eeq
where $b$ is the Bartlett correction constant and $(1-bn^{-1})$ is the Bartlett correction factor which depend the moments of $g(X,\theta_0)$. It follows from (\ref{2.36}) and  (\ref{2.37}) that
\beq
pr\{l(\theta_0)[1-bn^{-1}+O_p(n^{-3/2})]\leq c\}=pr(\chi^2_q \leq c)+O(n^{-2}). \label{2.40}
\eeq
Let $l_B(\theta)=(1-bn^{-1})l(\theta)$ be the Bartlett corrected empirical log-likelihood ratio, and denote by ${\mathcal C}'_{1-\alpha}$ the Bartlett corrected empirical likelihood confidence region for $\theta_0$.  Then,
\beq
{\mathcal C}'_{1-\alpha}=\{\theta: \theta \in \Theta_n \mbox{ and } l_B(\theta) \leq c\}.   \label{2.50}
\eeq
Equation (\ref{2.40}) implies that
\beq
pr(\theta_0 \in {\mathcal C}_{1-\alpha}')=P[l(\theta_0)(1-bn^{-1})\leq c]=pr(\chi^2_p \leq c)+O(n^{-2}). \label{2.52}
\eeq
A more detailed reviewed of the Bartlett correction is given the supplemental material.

\subsection{Composite similarity mapping}

The mismatch between the OEL domain $\Theta_n$ and the parameter space $\mathbb{R}^p$ is a main cause of the poor accuracy of the OEL confidence regions (Tsao, 2013). To solve the mismatch problem, we expand $\Theta_n$ to $\mathbb{R}^p$ through a {composite similarity mapping} $h^C_n: \Theta_n \rightarrow \mathbb{R}^p$  (Tsao and Wu, 2013). Under the three conditions, there exists a $\surd{n}$-consistent maximum empirical likelihood estimator $\tilde{\theta}$ for $\theta_0$ (see Appendix). 
Using OEL $l(\theta)$ and $\tilde{\theta}$, we define $h^C_n$ as
\begin{equation}
h_n^C(\theta)=\tilde{\theta}+\gamma(n,l(\theta))(\theta-\tilde{\theta}) \hspace{0.2in} \mbox{for \; $\theta \in \Theta_n$}, \label{3.1}
\end{equation}
where function $\gamma(n,l(\theta))$ is the {\em expansion factor} given by 
\begin{equation}
\gamma(n,l(\theta))=1+\frac{l(\theta)}{2n}. \label{3.2}
\end{equation}
To see how $h^C_n$ maps $\Theta_n$ to $\mathbb{R}^p$, define the level-$\tau$ contour of the OEL $l(\theta)$ as,
\beq
c(\tau)=\{\theta: \theta \in \Theta_n \hspace{0.07in} \mbox{and} \hspace{0.07in} l(\theta)=\tau\}, \label{3.10}
\eeq
where $\tau\geq \tilde{\tau}=l(\tilde{\theta})\geq 0$. For the just-determined case, $\tilde{\theta}$ is the solution of 
$\sum_{i=1}^ng(X_i,\theta)=0$, thus $R(\tilde{\theta})=1$ and $\tilde{\tau}=l(\tilde{\theta})=0$.  The contours form a partition of the OEL domain,
\beq
\Theta_n = \bigcup_{\tau\in[\tilde{\tau},+\infty)} c(\tau). \label{3.11}
\eeq 
Under the condition (which we will refer to as {\em condition 4}) that each OEL contour is the boundary of a connected region and the OEL contours are nested, (\ref{3.11}) implies that $c(\tilde{\tau})=\{\tilde{\theta}\}$ is the centre of $\Theta_n$. The value of $\tau$ measures the outwardness of a $c(\tau)$ with respect to the centre; the larger the $\tau$ value, the more outward $c(\tau)$ is. 
Theorem \ref{thm1} below gives the key properties of $h^C_n$.

\begin{theorem} \label{thm1}
Under conditions 1, 2 and 3, $h^C_n$ defined by (\ref{3.1}) and (\ref{3.2}) satisfies: 
\\ \hspace*{0.3in} (i) $h^C_n$ has a unique fixed point at $\tilde{\theta}$;
\\ \hspace*{0.3in} (ii) it is a similarity transformation for each individual OEL contour;
\\ \hspace*{0.3in} (iii) it is a surjection from $\Theta_n$ to $\mathbb{R}^p$.
\end{theorem}

Because of ($ii$), we call $h_n^C$ the composite similarity mapping as it may be viewed as a continuous sequence of similarity mappings from $\mathbb{R}^p$ to $\mathbb{R}^p$ indexed by $\tau\in [\tilde{\tau}, +\infty)$. 
The ``$\tau$-th'' mapping has expansion factor $\gamma(n,l(\theta))=\gamma(n,\tau)$ and is used exclusively to map the ``$\tau$-th'' OEL contour $c(\tau)$. Since $\gamma(n,\tau)$ is an increasing function of $\tau$, contours farther away from the centre are expanded more so that images of the contours fill up the entire $\mathbb{R}^p$. But regardless of the amount expanded, an OEL contour and its image are identical in shape; Figure 1 illustrates this with OEL contours for parameters of a regression model and their expanded images.

The proof of Theorem 1 is given in the supplementary material. A remark following the proof shows that if we are to add condition 4 to Theorem \ref{thm1}, then ($iii$) can be strengthened to ($iii'$) $h^C_n$ is a bijection from $\Theta_n$ to $\mathbb{R}^p$. It is not clear how we may verify condition 4 through $g(X,\theta)$. This is why we have kept it separate from the three conditions identified in the preliminaries. Nevertheless, we have not encountered any examples where condition 4 is violated. 

\subsection{Extended empirical likelihood on full parameter space}

By Theorem \ref{thm1}, $h_n^C: \Theta_n\rightarrow \mathbb{R}^d$ is surjective. Thus, for any $\theta \in \mathbb{R}^p$,  $s(\theta)=\{\theta': h_n^C(\theta')=\theta\}$ is non-empty. When $h_n^C$ is not injective, $s(\theta)$ may contain more than one point and $h_n^C$ does not have an inverse. Hence, we define a generalized inverse $h_n^{-C}: \mathbb{R}^p \rightarrow  \Theta_n$ as follows,
\beq
h_n^{-C}(\theta)=\argmin_{\theta '\in s(\theta)} \{\|\theta'-\theta\|\}. \label{3.29}
\eeq
The extended empirical log-likelihood ratio EEL $l^*(\theta)$ under $h_n^{-C}$ is then
\beq
l^*(\theta)=l(h_n^{-C}(\theta)) \hspace{0.2in} \mbox{for $\theta \in \mathbb{R}^p$},  \label{3.30}
\eeq
which is well-defined throughout $\mathbb{R}^p$. We now give the properties of the point $\theta_0'$ satisfying
\beq
h_n^{-C}(\theta_0)=\theta_0', \label{3.20}
\eeq 
and the asymptotic distribution of $l^*(\theta_0)=l(h_n^{-C}(\theta_0))=l(\theta_0')$. For convenience, we use $[\tilde{\theta}, \theta_0]$ to denote the line segment that connects $\tilde{\theta}$ and $\theta_0$. We have

\begin{lemma}  \label{lem1} 
Under conditions 1, 2 and 3, point $\theta_0'$ defined by equation (\ref{3.20}) satisfies \\[0.1in]
\hspace*{1.2in} ($i$) $\theta_0' \in [\tilde{\theta}, \theta_0]$ \hspace{0.05in} and \hspace{0.05in} ($ii$) $\theta'_0-\theta_0 = O_p(n^{-3/2})$. \end{lemma}

\begin{theorem}  \label{thm2}
Under conditions 1, 2 and 3, the EEL $l^*(\theta)$ defined by (\ref{3.30})  satisfies
\beq
l^*(\theta_0) {\longrightarrow} \chi^2_q  \label{3-40}
\eeq 
\hspace*{1in} in distribution as $n\rightarrow +\infty$.
\end{theorem}

Proofs of Lemma 1 and Theorem \ref{thm2} are sketched in the Appendix. Detailed proofs are given in the supplementary material. A key element in the proof for Theorem \ref{thm2} is the following simple relationship between the OEL $l(\theta)$ and the EEL $l^*(\theta)$:
\beq
l^*(\theta_0)=l(h_n^{-C}(\theta_0))= l(\theta'_0)=l(\theta_0+(\theta_0'-\theta_0)).  \label{3.35}
\eeq
This and the fact that $\|\theta_0'-\theta_0\|$ is asymptotically very small imply that $l^*(\theta_0)=l(\theta_0)+o_p(1)$. Relation (\ref{3.35}) is also the key in the derivation of a second order EEL in the next section.

\subsection{Second order extended empirical likelihood}

The BEL of DiCiccio, Hall and Romano (1991) has the second order accuracy. Theorem \ref{thm3} shows that for the just-determined case the EEL can also attain the second order accuracy.

\begin{theorem} \label{thm3} Assume conditions 1, 2 and 3 hold. For the just-determined case where $p=q$, let $l^*_2(\theta)$ be the EEL under the composite similarity mapping (\ref{3.1}) with expansion factor
\begin{equation}
\gamma_2(n,l(\theta))=1+\frac{b }{2n}[l(\theta)]^{\delta(n)}, \label{4.15}
\end{equation}
where $\delta(n)=O(n^{-1/2})$ and $b$ is the Bartlett correction constant in (\ref{2.37}) and (\ref{2.40}). Then 
\beq
l^*_2(\theta_0)=l(\theta_0)[1-bn^{-1}+O_p(n^{-3/2})]. \label{4.20}
\eeq
\end{theorem}
Proof of Theorem \ref{thm3} is given in the supplementary material. Comparing (\ref{4.20}) with (\ref{2.40}), we see that $l^*_2(\theta)$ is equivalent to the BEL $l_B(\theta)$. Hence, we call it the {\em second order EEL}. Correspondingly, we call $l^*(\theta)$ defined by the $\gamma(n, l(\theta))$ in (\ref{3.2}) the {\em first order EEL}. The utility of the $\delta(n)$ in $\gamma_2(n,l(\theta))$ is to control the speed of domain expansion which ensures $l^*_2(\theta)$ behaves asymptotically like $l_B(\theta)$. For convenience, in our numerical comparison we set $\delta(n)=n^{-1/2}$. 

We noted after Theorem \ref{thm2} that $l^*(\theta_0)=l(\theta_0)+o_p(1)$. An even stronger connection between $l^*(\theta_0)$ and $l(\theta_0)$ is given by Corollary \ref{cor1} below. This result helps to explain the remarkable numerical accuracy of confidence regions based on the first order EEL $l^*(\theta)$ in the next section.

\begin{corollary} \label{cor1}  Under conditions 1, 2 and 3, EEL $l^*(\theta)$ for the just-determined case satisfies,
\beq
l^*(\theta_0)=l(\theta_0)[1-l(\theta_0)n^{-1}+O_p(n^{-3/2})]. \label{4.30}
\eeq
\end{corollary}

\section{Numerical examples}  

We compare the first order EEL $l^*(\theta)$ with the first order OEL and the second order BEL through a small simulation study. A more comprehensive comparison is given in the supplementary material. 
Table \ref{table1} contains simulated coverage probabilities for $\boldsymbol{\beta}$ of linear model
\[ y=\mathbf{x}^T\boldsymbol{\beta} +\varepsilon,\]
where $\varepsilon\sim N(0,1)$. We consider two models: Model 1 given by $\mathbf{x}=(1,x_1)^T$ and $\boldsymbol{\beta}=(1,2)^T$ and Model 2 given by $\mathbf{x}=(1,x_1,x_2)^T$ and $\boldsymbol{\beta}=(1,2,3)^T$. For the simulation, values of $x_1$ are randomly generated from a uniform distribution on $[0,30]$ and that of $x_2$ are randomly generated from a uniform distribution on $[20,50]$. The EEL methods are defined by the composite similarity mapping centred on $\tilde{\theta}=\hat{\boldsymbol{\beta}}$, the least-squares estimate of $\boldsymbol{\beta}$.

{\em Model 1 based comparison:} The EEL is consistently more accurate than the OEL for all combinations of sample size and confidence level. In particular, for small to moderate sample sizes ($n\leq 30$) it is substantially more accurate than the OEL. The EEL is also more accurate than the BEL for small to moderate sample sizes. Remarkably, even for large sample sizes ($n>30$), it remains more accurate than the second order BEL. This surprising observation may be partially explained by Corollary 1 where the EEL is seen as having a Bartlett correction type of expansion. See the supplementary material for more examples and further discussion. 

{\em Model 2 based comparison:} The parameter vector of Model 2 has dimension $p=3$ whereas that of Model 1 has $p=2$. This difference allows us to assess the impact of dimension $p$. When $p$ increases from 2 to 3, the coverage probability of the EEL is the least affected. For small to moderate sample sizes, that of the OEL and BEL deteriorated a lot. This is due to the mismatch problem which has bigger impact on the OEL and BEL in higher dimensions. The EEL is not affected by the mismatch, thus it held up much better. In particular, the 99\% EEL confidence region is the most reliable and is accurate for all combinations of $n$ and $p$.

We conclude by briefly commenting on the computation of EEL $l^*(\theta)$. Suppose $h_n^C$ is also injective. Since $l^*(\theta)=l(\theta')$, we compute $l^*(\theta)$ by finding the $\theta'$ satisfying $h_n^{C}(\theta')=\theta$ first and then compute $l(\theta')$. We may find $\theta'$ by computing the root for the multivariate function $f(\theta')=h_n^{C}(\theta')-\theta$. But it is more efficient to reformulate this function as a {\em univariate} function by using the fact that $\theta' \in [\tilde{\theta}, \theta]$ (see Theorem \ref{thm1} and its proof). When $h_n^C$ is not necessarily injective, we find one $\theta'$ satisfying $h_n^{C}(\theta')=\theta$ first (call it $\theta_1'$). Then, look for another satisfying $h_n^{C}(\theta')=\theta$ in the interval $[\theta_1', \theta]$, and iterate this process until no new solutions can be found. The last of these (call it $\theta_l'$) is the solution closest to $\theta$ and hence $l^*(\theta)=l(\theta_l')$.

\begin{table}
\def~{\hphantom{0}}
\tbl{Coverage probabilities of OEL, EEL and BEL confidence regions}{%
\begin{tabular}{llccccccccc} \\
&    & \multicolumn{3}{c}{90\% level} & \multicolumn{3}{c}{95\% level} & \multicolumn{3}{c}{99\% level}\\[1pt]
& $n$ & OEL & EEL & BEL & OEL & EEL & BEL & OEL & EEL & BEL \\[5pt]
\textbf{Model 1} & 10 & 66.9 & 80.0 & 76.3 & 73.4 & 88.5 & 80.9 & 81.5 & 98.4 & 87.5\\
                 & 20 & 79.7 & 85.6 & 85.1 & 86.5 & 92.5 & 90.8 & 94.3 & 98.5 & 96.6\\
                 & 30 & 84.3 & 87.8 & 87.2 & 90.1 & 93.9 & 92.6 & 96.5 & 98.6 & 97.5\\
                 & 50 & 86.7 & 88.8 & 88.5 & 92.6 & 94.3 & 93.7 & 97.7 & 98.9 & 98.2\\
                & 100 & 88.8 & 89.8 & 89.6 & 94.0 & 94.8 & 94.5 & 98.4 & 99.0 & 98.6\\[3pt]
\textbf{Model 2} & 10 & 47.3 & 75.1 & 58.6 & 54.1 & 87.2 & 64.8 & 65.1 & 97.7 & 74.2\\
                 & 20 & 69.9 & 81.2 & 77.6 & 77.3 & 89.7 & 84.2 & 88.0 & 97.8 & 92.3\\
                 & 30 & 76.8 & 84.3 & 83.0 & 84.4 & 91.1 & 88.8 & 92.9 & 98.1 & 95.5\\
                 & 50 & 83.5 & 87.2 & 86.8 & 89.8 & 93.1 & 92.0 & 96.3 & 98.5 & 97.6\\
                & 100 & 87.4 & 89.1 & 88.8 & 93.0 & 94.4 & 94.0 & 98.4 & 99.0 & 98.6\\
\end{tabular}}
\label{table1}
\begin{tabnote}
Each entry in the table is a simulated coverage probability for $\boldsymbol{\beta}$ based on 10,000 random samples of size $n$ indicated in column 2 from the linear model indicated in column 1.
\end{tabnote}
\end{table}


\section{Discussion}

The impressive accuracy of the first order EEL can also be seen through the examples in the supplementary material. We recommend it for practical applications due to its simplicity and superior accuracy. Although the focus of this paper is on EEL for parameters defined by estimating equations, main techniques employed in the proofs may be applied to handle parameters in other settings. In general, an EEL may be derived so long as a $\surd{n}$-consistent maximum empirical likelihood estimator $\tilde{\theta}$ is available. If the OEL contours are nested, then the EEL retains not only all asymptotic properties of the OEL but also the geometric characteristics of its contours.  Finally, we have only considered the case where the full parameter space $\Theta$ is $\mathbb{R}^p$. The case where $\Theta$ is a known subset of $\mathbb{R}^p$ may be handled by finding the EEL on $\mathbb{R}^p$ first, and then redefining it as positive infinity for $\theta\notin \Theta$ while keeping it unchanged for $\theta\in \Theta$.

\section*{Acknowledgement}
The first author is supported by a grant from the National Science and Engineering Research Council of Canada.

\section*{Supplementary material}
\label{SM}
Supplementary material available 
online includes detailed proofs of all Lemmas and Theorems (Part I) and a more comprehensive numerical comparison (Part II).

\appendix  

\appendixone
\section*{Appendix}

We identify two assumptions which are used implicitly in the proofs. We also sketch the proofs for Lemma 1 and Theorem 2. Detailed proofs are provided in the supplementary material.

Under conditions 1, 2 and 3, we assume without loss of generality that (a) the OEL domain $\Theta_n$ is an open set containing $\theta_0$ and (b) there exists a $\surd{n}$-consistent maximum empirical likelihood estimator $\tilde{\theta}$ for $\theta_0$. To see (a), by condition 1 and Lemma 11.1 in Owen (2001), with probability tending to 1 that the convex hull of the $g(X_i,\theta_0)$ contains 0. Hence, we may assume for sufficiently large $n$ that $\Theta_n$ contains $\theta_0$. To see $\Theta_n$ is also open, suppose $\theta\in \Theta_n$. Then, the convex hull of the $g(X_i,\theta)$ contains $0$ in its interior. By condition 2, $g(X_i,\theta)$ is continuous in $\theta$ which implies that a small change in $\theta$ will result in only a small change in the convex hull. Thus, there exists a small neighbourhood of $\theta$ such that for any $\theta'$ in that neighbourhood the convex hull of the $g(X_i,\theta')$ also contains $0$. Hence, this neighbourhood is inside $\Theta_n$ and $\Theta_n$ is open. To see (b), we refer to Lemma 1 and Theorem 1 in Qin and Lawless (1994) which give, respectively, the existence and $\surd{n}$-consistency of the maximum empirical likelihood estimator.

\vspace{0.1in}

\begin{proof}[of Lemma 1]

Differentiating both sides of equation (\ref{2.15}) with respect to $\theta$, we obtain $J(\theta_0)=[\partial l(\theta)/\partial \theta]_{\theta=\theta_0}=O_p(n^{1/2})$.
For $\theta$ values in a small neighbourhood of $\theta_0$, $\{\theta: \|\theta-\theta_0\|\leq \kappa n^{-1/2}\}$, where $\kappa$ is a positive constant, Taylor expansion of $l(\theta)$ gives
\begin{equation}
 l(\theta)=l(\theta_0+(\theta-\theta_0))=l(\theta_0)+J(\theta_0)(\theta-\theta_0) + O_p(1). \label{a-1}
 \end{equation}
Since $J(\theta_0)=O_p(n^{1/2})$ and $l(\theta_0)=O_p(1)$, (\ref{a-1}) implies that $l(\theta)=O_p(1)$. Also, $\gamma(n,l(\theta))\geq 1$ and
\beq  \theta_0-\tilde{\theta}=\gamma(n,l(\theta_0'))(\theta_0'-\tilde{\theta}), \label{a-5} \eeq
thus $\theta'_0$ is on the ray originating from $ \tilde{\theta}$ through $\theta_0$ and $\|\theta_0-\tilde{\theta}\|\geq\|\theta_0'-\tilde{\theta}\|$. Hence, $\theta_0'\in [\tilde{\theta}, \theta_0]$. This and the $\surd{n}$-consistency of $\tilde{\theta}$ imply that $\theta_0'-\theta_0 =O_p(n^{-1/2})$. It follow that $l(\theta_0')=O_p(1)$ and
\begin{equation}
\gamma(n,l(\theta_0'))=1+\frac{l(\theta_0')}{2n}=1+O_p(n^{-1}). \nonumber
\end{equation}
This and (\ref{a-5}) then imply $\theta_0'-\theta_0=O_p(n^{-3/2})$.
\end{proof}

\begin{proof}[of Theorem 2]
By Lemma \ref{lem1} ($ii$), $\theta_0'-\theta_0=O_p(n^{-3/2})$. Taylor expansion of $l^*(\theta_0)$ gives
\begin{equation}
 l^*(\theta_0)=l(\theta_0')=l(\theta_0+(\theta_0'-\theta_0))=l(\theta_0)+J(\theta_0)(\theta_0'-\theta_0) + O_p(n^{-1}). \label{a-2}
 \end{equation}
Since $J(\theta_0)=O_p(n^{1/2})$, (\ref{a-2}) implies that 
$l^*(\theta_0)=l(\theta_0)+ O_p(n^{-1})$. Hence, the EEL $l^*(\theta_0)$ has the same limiting $\chi^2_q$ distribution as the OEL $l(\theta_0)$. 
\end{proof}


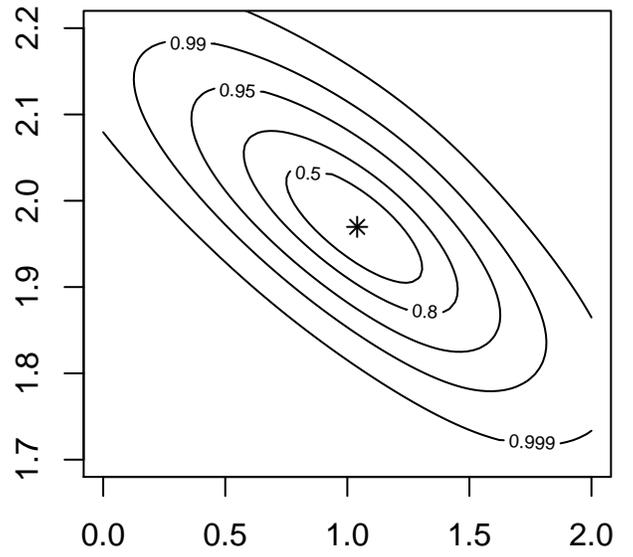
\end{document}